# PRODUCTS OF BLOCK TOEPLITZ OPERATORS

CAIXING GU AND DECHAO ZHENG

ABSTRACT. In this paper we characterize when the product of two block Toeplitz operators is a compact perturbation of a block Toeplitz operator on the Hardy space of the open unit disk. Necessary and sufficient conditions are given for the commutator of two block Toeplitz operators to be compact.

## 1. INTRODUCTION

Let $D$ be the open unit disk in the complex plane and $\partial D$ the unit circle. Let $d\sigma(w)$ be the normalized Lebesgue measure on the unit circle. We denote by $L^2(C^n)$ ($L^2$ for $n=1$) the space of $C^n$-valued Lebesgue square integrable functions on the unit circle. The Hardy space $H^2(C^n)$ ($H^2$ for $n=1$) is the closed linear span of $C^n$-valued analytic polynomials. We observe that $L^2(C^n) = L^2 \otimes C^n$ and $H^2(C^n) = H^2 \otimes C^n$, where $\otimes$ denotes the Hilbert space tensor product. Let $M_{n \times n}$ be the set of $n \times n$ complex matrices. $L^\infty_{n \times n}$ denotes the space of $M_{n \times n}$-valued essentially bounded Lebesgue measurable functions on the unit circle and $H^\infty_{n \times n}$ denotes the space of $M_{n \times n}$-valued essentially bounded analytic functions in the disk.

Let $P$ be the projection of $L^2(C^n)$ onto $H^2(C^n)$. For $F \in L^\infty_{n \times n}$, the block Toeplitz operator $T_F : H^2(C^n) \to H^2(C^n)$ with symbol $F$ is defined by the rule $T_F h = P(Fh)$. The Hankel operator $H_F : H^2(C^n) \to L^2(C^n) \ominus H^2(C^n)$ with symbol $F$ is defined by $H_F h = (I-P)(Fh)$. The block Toeplitz operator $T_F$ has the following matrix representation:

$$\begin{bmatrix} A_0 & A_{-1} & A_{-2} & \cdots \\ A_1 & A_0 & A_{-1} & \cdots \\ A_2 & A_1 & A_0 & \cdots \\ \vdots & \vdots & \vdots & \vdots \end{bmatrix}$$

where $A_i$ belongs to $M_{n \times n}$. The word "block" refers to the fact that in the above matrix representation the entries are not scalars but linear transformations on $C^n$. In this paper the word "block" will often be omitted. For more details on the block Toeplitz operators and Hankel operators, see [7], [9] and [3].

---

1991 *Mathematics Subject Classification.* Primary 47B35.

Gu's research is supported in part by NSF Grant DMS-9022140 during residence at MSRI. Zheng is supported in part by a NSF grant.

If we set $H^2(C^n) = H^2 \oplus ... \oplus H^2$, then we see that the Toeplitz operator $T_F$ has the form

$$T_F = \begin{bmatrix} T_{f_{11}} & T_{f_{12}} & \cdots & T_{f_{1n}} \\ T_{f_{21}} & T_{f_{22}} & \cdots & T_{f_{2n}} \\ \vdots & \vdots & \vdots & \vdots \\ T_{f_{n1}} & T_{f_{n2}} & \cdots & T_{f_{nn}} \end{bmatrix}$$

and the Hankel operator $H_F$ has the form

$$H_F = \begin{bmatrix} H_{f_{11}} & H_{f_{12}} & \cdots & H_{f_{1n}} \\ H_{f_{21}} & H_{f_{22}} & \cdots & H_{f_{2n}} \\ \vdots & \vdots & \vdots & \vdots \\ H_{f_{n1}} & H_{f_{n2}} & \cdots & H_{f_{nn}} \end{bmatrix},$$

where

$$F = \begin{bmatrix} f_{11} & f_{12} & \cdots & f_{1n} \\ f_{21} & f_{22} & \cdots & f_{2n} \\ \vdots & \vdots & \vdots & \vdots \\ f_{n1} & f_{n2} & \cdots & f_{nn} \end{bmatrix}.$$

Thus as in the scalar case, the block Toeplitz operators and Hankel operators are connected by the following important relation:

$$T_{FG} - T_F T_G = H_{F^*}^* H_G.$$

The map $\xi \colon F \to T_F$, which is called the Toeplitz quantization, carries $L_{n \times n}^\infty$ into the $C^*$−algebra of bounded operators on $H^2(C^n)$. It is a contractive *-linear mapping [7]. However it is not multiplicative in general. Indeed Brown and Halmos [4] showed that for scalar functions $f$ and $g$, $T_f T_g = T_{fg}$ if and only either $f^*$ or $g$ is in $H^\infty$. It is not difficult to see that in the matrix case $T_F T_G = T_{FG}$ still holds if either $F^*$ or $G$ is in $H_{n \times n}^\infty$. But the converse is not valid in the matrix case. We will characterize $F$ and $G$ such that $T_F T_G = T_{FG}$; see Theorem 6 below for details. On the other hand, Douglas [7] showed that $\xi$ is actually a cross section for a *-homomorphism from the Toeplitz algebra, the $C^*$−algebra generated by all bounded Toeplitz operators on $H^2(C^n)$, onto $L_{n \times n}^\infty$. So modulo the commutator ideal of the Toeplitz algebra, $\xi$ is multiplicative.

The main question to be considered in this paper is when the product of two Toeplitz operators is a compact perturbation of a Toeplitz operator. This problem is connected with the spectral theory of Toeplitz operators; see [7], [9] and [3]. It follows from a theorem of Douglas [7] that $T_F T_G$ can be a compact perturbation of a Toeplitz operator only when it is a compact perturbation of $T_{FG}$. Thus it suffices to study when the semi-commutator $T_{FG} - T_F T_G$ is compact. When $F = f$ and $G = g$ are scalar functions, the problem was solved by Axler, Chang, and Sarason [2] and Volberg [14]. Their beautiful result is that $T_{fg} - T_f T_g$ is compact if and only if $H^\infty[\bar{f}] \bigcap H^\infty[g] \subset H^\infty + C(\partial D)$; here $H^\infty[g]$ denotes the closed subalgebra of $L^\infty$ generated by $H^\infty$ and $g$.

Recently, Zheng [15] proved that $T_{fg} - T_f T_g$ is compact if and only if

$$\lim_{|z| \to 1} \|H_{\bar{f}} k_z\|_2 \|H_g k_z\|_2 = 0;$$

here $k_z$ denotes the normalized reproducing kernel in $H^2$ for point evaluation at $z$. If we write $f = f_+ + f_-$ for each $f \in L^\infty$ where $f_+$ and $\overline{f_-}$ are in $H^2$, then the above condition is equivalent to

$$\lim_{|z| \to 1} |\overline{f_+} - \overline{f_+}(z)|^2(z)|g_- - g_-(z)|^2(z) = 0,$$

where $h(z)$ denotes the harmonic extension of $h$ at $z \in D$ for $h \in L^1$, via the Poisson integral

$$h(z) = \int_{\partial D} h(w) \frac{(1 - |z|^2)}{|1 - w\bar{z}|^2} d\sigma(w).$$

For the block Toeplitz operators, we will show that $T_{FG} - T_F T_G$ is compact if and only if

$$\lim_{|z| \to 1} \|[|(F_+)^* - (F_+)^*(z)|^2(z)]^{1/2}[|G_- - G_-(z)|^2(z)]^{1/2}\| = 0,$$

where we write $F = (f_{ij})_{n \times n}$ as $F = F_+ + F_-$ with $F_+ = ((f_{ij})_+)_{n \times n}$ and $F_- = ((f_{ij})_-)_{n \times n}$, and

$$H(z) = (h_{ij}(z))_{n \times n}$$

if $H = (h_{ij})_{n \times n}$. For a matrix $A$, we define $|A|^2 = AA^*$. Several other equivalent conditions, in particular a condition in the spirit of the result of Axler, Chang, and Sarason [2] and Volberg [14], will be given.

In [11] Gorkin and Zheng characterized when the commutator $T_f T_g - T_g T_f$ is compact for scalar functions $f$ and $g$. In this paper, by considering block Toeplitz operators, we will give a unified approach for the study of compactness of both semi-commutators and commutators. Namely, by a theorem of Douglas [7], we have that the commutator $T_F T_G - T_G T_F$ of block Toeplitz operators $T_F$ and $T_G$ is compact if and only if $FG = GF$ and the semi-commutator $T_{BC} - T_B T_C$ is compact, where

$$B = \begin{bmatrix} F & -G \\ 0 & 0 \end{bmatrix}, \quad C = \begin{bmatrix} G & 0 \\ F & 0 \end{bmatrix};$$

see Theorem 7 below for details. Thus we will show the commutator $T_F T_G - T_G T_F$ is compact if and only if $FG = GF$ and

$$\lim_{|z| \to 1} \left\| \begin{bmatrix} |(F_+)^* - (F_+)^*(z)|^2(z) & -(((F_+)^* - (F_+)^*(z))(G_+ - G_+(z)))(z) \\ -(((G_+)^* - (G_+)^*(z))(F_+ - F_+(z)))(z) & |(G_+)^* - (G_+)^*(z)|^2(z) \end{bmatrix}^{1/2} \times \right.$$

$$\left. \begin{bmatrix} |G_- - G_-(z)|^2(z) & ((G_- - G_-(z))((F_-)^* - (F_-)^*(z)))(z) \\ ((F_- - F_-(z))((G_-)^* - (G_-)^*(z)))(z) & |F_- - F_-(z)|^2(z) \end{bmatrix}^{1/2} \right\| = 0.$$

## 2. A necessary condition for compactness

In this section we will obtain a necessary condition for the compactness of the semi-commutator $T_{FG} - T_F T_G$. This will also motivate a necessary and sufficient condition for $T_{FG} - T_F T_G$ to be zero. The question when $T_{FG} - T_F T_G$ is zero will be discussed in the next two sections.

First we introduce an antiunitary operator $V$ on $L^2$ by defining $(Vh)(w) = \overline{wh(w)}$. The operator enjoys many nice properties such as $V^{-1}(I-P)V = P$ and $V = V^{-1}$. These properties leads easily to the relation $V^{-1}H_f V = H_f^*$.

Let $x$ and $y$ be two vectors in $L^2$. $x \otimes y$ is the operator of rank one defined by

$$(x \otimes y)(f) = <f, y> x.$$

Observe that the norm of the operator $x \otimes y$ is $\|x\|_2 \|y\|_2$.

For $z$ in $D$, let $k_z$ be the normalized reproducing kernel $\frac{(1-|z|)^{1/2}}{(1-\overline{z}w)}$ for point evaluation at $z$, and $\phi_z$ the Möbius map on the unit disk,

$$\phi_z(w) = \frac{z-w}{1-\overline{z}w}.$$

$\phi_z$ can also be viewed as a function on the unit circle. Let $\Phi_z$ denote the function $diag\{\phi_z, \cdots, \phi_z\} \in H_{n \times n}^\infty$. The product $T_{\phi_z} T_{\overline{\phi_z}}$ is the orthogonal projection onto $H^2 \ominus \{k_z\}$. Thus $1 - T_{\phi_z} T_{\overline{\phi_z}}$ is the operator $k_z \otimes k_z$ of rank one. This leads to the following lemma.

**Lemma 1.** *Let $F = (f_{ij})_{n \times n}$ and $G = (g_{ij})_{n \times n}$ be in $L_{n \times n}^\infty$. Then the operator $H_F^* H_G - T_{\Phi_z}^* H_F^* H_G T_{\Phi_z}$ is anti-unitary equivalent to*

$$(\sum_{j=1}^n H_{f_{ji}} k_z \otimes H_{g_{jk}} k_z)_{n \times n}.$$

*Proof.* Let $F = (f_{ij})_{n \times n}$ and $G = (g_{ij})_{n \times n}$. Then it is easy to check that

$$H_F^* H_G = (\sum_{j=1}^n H_{f_{ji}}^* H_{g_{jk}})_{n \times n}$$

and

$$T_{\Phi_z}^* H_F^* H_G T_{\Phi_z} = (\sum_{j=1}^n T_{\phi_z}^* H_{f_{ji}}^* H_{g_{jk}} T_{\phi_z})_{n \times n}.$$

So the difference $H_F^* H_G - T_{\Phi_z}^* H_F^* H_G T_{\Phi_z}$ is

$$(\sum_{j=1}^n [H_{f_{ji}}^* H_{g_{jk}} - T_{\phi_z}^* H_{f_{ji}}^* H_{g_{jk}} T_{\phi_z}])_{n \times n}.$$

Hence it is sufficient to check that

$$\sum_{j=1}^{n} [H^*_{f_{ji}} H_{g_{jk}} - T^*_{\phi_z} H^*_{f_{ji}} H_{g_{jk}} T_{\phi_z}]$$

is unitary equivalent to

$$\sum_{j=1}^{n} H_{f_{ji}} k_z \otimes H_{g_{jk}} k_z.$$

Applying $V$ and $V^{-1}$ to the difference

$$\sum_{j=1}^{n} [H^*_{f_{ji}} H_{g_{jk}} - T^*_{\phi_z} H^*_{f_{ji}} H_{g_{jk}} T_{\phi_z}]$$

we have

$$V^{-1} \sum_{j=1}^{n} [H^*_{f_{ji}} H_{g_{jk}} - T^*_{\phi_z} H^*_{f_{ji}} H_{g_{jk}} T_{\phi_z}] V$$

$$= V^{-1} \sum_{j=1}^{n} [H^*_{f_{ji}} H_{g_{jk}} - H^*_{f_{ji}\phi_z} H_{g_{jk}\phi_z}] V.$$

Because $V^{-1} H_f V = H^*_f$, the above equality simplifies to

$$= V^{-1} \sum_{j=1}^{n} [H^*_{f_{ji}} H_{g_{jk}} - H^*_{f_{ji}\phi_z} H_{g_{jk}\phi_z}] V$$

$$= \sum_{j=1}^{n} [H_{f_{ji}} H^*_{g_{jk}} - H_{f_{ji}\phi_z} H^*_{g_{jk}\phi_z}]$$

$$= \sum_{j=1}^{n} H_{f_{ji}} (1 - T_{\phi_z} T^*_{\phi_z}) H^*_{g_{jk}}.$$

Since $1 - T_{\phi_z} T^*_{\phi_z} = k_z \otimes k_z$, the right hand side of the above equation is

$$\sum_{j=1}^{n} H_{f_{ji}} k_z \otimes H_{g_{jk}} k_z$$

for all $z$ in $D$. This completes the proof of the lemma. ∎

Let $trace$ be the trace on the trace class of operators on $H^2$ and $tr$ denote the trace on the n by n matrices.

**Lemma 2.** *Let $F$ and $G$ be in $L^\infty_{n \times n}$. Let $T = H^*_F H_G - T^*_{\Phi_z} H^*_F H_G T_{\Phi_z}$. Then*

$$trace\{T^*T\} = tr[|F_- - F_-(z)|^2(z)|G_- - G_-(z)|^2(z)].$$

*Proof.* Let $F = (f_{ij})_{n \times n}$ and $G = (g_{ij})_{n \times n}$. By Lemma 1, $H_F^* H_G - T_{\Phi_z}^* H_F^* H_G T_{\Phi_z}$ is anti-unitary equivalent to

$$T_1 = (\sum_{j=1}^n H_{f_{ji}} k_z \otimes H_{g_{jk}} k_z)_{n \times n}.$$

So we need to computer the trace of the operator $T_1^* T_1$. It is easy to see that

$$T_1^* T_1 = (\sum_{j_1} (\sum_j H_{f_{jj_1}} k_z \otimes H_{g_{jk}} k_z)^* (\sum_\mu H_{f_{\mu j_1}} k_z \otimes H_{g_{\mu l}} k_z))$$

$$= (\sum_{j_1} \sum_j \sum_\mu < H_{f_{\mu j_1}} k_z, H_{f_{jj_1}} k_z > H_{g_{jk}} k_z \otimes H_{g_{\mu l}} k_z),$$

where the second equality above follows from the fact that

$$(x_1 \otimes y_1)^* (x_2 \otimes y_2) = < x_2, x_1 > y_1 \otimes y_2.$$

Thus

$$trace\{T_1^* T_1\} = \sum_l \sum_{j_1} \sum_j \sum_\mu < H_{f_{\mu j_1}} k_z, H_{f_{jj_1}} k_z > < H_{g_{ji}} k_z, H_{g_{\mu l}} k_z >.$$

If we write $f_{ij} = (f_{ij})_+ + (f_{ij})_-$ and $g_{ij} = (g_{ij})_+ + (g_{ij})_-$ for $(f_{ij})_+$, $(g_{ij})_+ \in H^2$ and $(f_{ij})_-$, $(g_{ij})_- \in \overline{zH^2}$, then $H_{f_{ij}} k_z = ((f_{ij})_- - (f_{ij})_-(z)) k_z$ and $H_{g_{ij}} k_z = ((g_{ij})_- - (g_{ij})_-(z)) k_z$. Therefore, by a change of the order of summations, we have

$$trace\{T_1^* T_1\} = \sum_l \sum_{j_1} \sum_j \sum_\mu < [(f_{\mu j_1})_- - (f_{\mu j_1})_-(z)] k_z, [(f_{jj_1})_- - (f_{jj_1})_-(z)] k_z >$$

$$\times < [(g_{jl})_- - (g_{jl})_-(z)] k_z, [(g_{\mu l})_- - (g_{\mu l})_-(z)] k_z >$$

$$= \sum_\mu \sum_j (\sum_{j_1} < [(f_{\mu j_1})_- - (f_{\mu j_1})_-(z)] k_z, [(f_{jj_1})_- - (f_{jj_1})_-(z)] k_z >) \times$$

$$(\sum_l < [(g_{jl})_- - (g_{jl})_-(z)] k_z, [(g_{\mu l})_- - (g_{\mu l})_-(z)] k_z >).$$

Note that

$$(\sum_{j_1} < [(f_{\mu j_1})_- - (f_{\mu j_1})_-(z)] k_z, [(f_{jj_1})_- - (f_{jj_1})_-(z)] k_z >)_{n \times n} = |F_- - F_-(z)|^2(z),$$

and similarly

$$(\sum_l < [(g_{jl})_- - (g_{jl})_-(z)] k_z, [(g_{\mu l})_- - (g_{\mu l})_-(z)] k_z >)_{n \times n} = |G_- - G_-(z)|^2(z).$$

Hence

$$trace\{T_1^* T_1\} = \sum_u \sum_j (|F_- - F_-(z)|^2(z))_{\mu j} (|G_- - G_-(z)|^2(z))_{j\mu}$$

$$= \sum_u [|F_- - F_-(z)|^2(z)) (|G_- - G_-(z)|^2(z))]_{\mu\mu}$$

$$= tr[|F_- - F_-(z)|^2(z)|G_- - G_-(z)|^2(z)].$$

Here we note that for a matrix $A$, $|A|^2 = AA^*$ and $(A)_{ij}$ denotes the $(i, j) - th$ entry of $A$. This completes the proof of this lemma. ∎

**Theorem 3.** *Let $F$ and $G$ be in $L^\infty_{n \times n}$. If $H^*_F H_G$ is compact, then*

$$\lim_{|z| \to 1} \|\{[|F_- - F_-(z)|^2(z)]^{1/2}[|G_- - G_-(z)|^2(z)]^{1/2}\}\| = 0.$$

*Proof.* Set $H^2(C^n) = H^2 \oplus ... \oplus H^2$. Let $(A_{ij}) := H^*_F H_G$ be the operator-valued $n \times n$ matrix representation of $H^*_F H_G$ with respect to the above decomposition of $H^2(C^n)$. Since $H^*_F H_G$ compact, each entry $A_{ij}$ of $H^*_F H_G$ is compact on $H^2$, by Lemma 2 [15], we have

$$\lim_{|z| \to 1} \|A_{ij} - T^*_{\phi_z} A_{ij} T_{\phi_z}\| = 0.$$

Hence

$$\lim_{|z| \to 1} \|H^*_F H_G - T^*_{\Phi_z} H^*_F H_G T_{\Phi_z}\| = 0.$$

By Lemma 1, $H^*_F H_G - T^*_{\Phi_z} H^*_F H_G T_{\Phi_z}$ is a finite rank operator. Therefore, we have

$$\lim_{|z| \to 1} trace\{(H^*_F H_G - T^*_{\Phi_z} H^*_F H_G T_{\Phi_z})^*(H^*_F H_G - T^*_{\Phi_z} H^*_F H_G T_{\Phi_z})\} = 0$$

because the norm of a finite rank positive operator is equivalent to its trace. By Lemma 2, we obtain

$$\lim_{|z| \to 1} tr[|F_- - F_-(z)|^2(z)|G_- - G_-(z)|^2(z)] = 0.$$

On the other hand,

$$tr[|F_- - F_-(z)|^2(z)|G_- - G_-(z)|^2(z)]$$

$$= tr\{[|F_- - F_-(z)|^2(z)]^{1/2}[|G_- - G_-(z)|^2(z)][|F_- - F_-(z)|^2(z)]^{1/2}\}$$

$$= tr\{[|F_- - F_-(z)|^2(z)]^{1/2}[|G_- - G_-(z)|^2(z)]^{1/2} \times$$

$$\{[|F_- - F_-(z)|^2(z)]^{1/2}[|G_- - G_-(z)|^2(z)]^{1/2}\}^*\}.$$

As is well-known, for all $n \times n$ matrices $A$,

$$tr A^* A \geq C\|A\|^2$$

for some constant $C > 0$. Hence we conclude

$$\lim_{|z| \to 1} \|\{[|F_- - F_-(z)|^2(z)]^{1/2}[|G_- - G_-(z)|^2(z)]^{1/2}\}\| = 0.$$

This completes the proof of the theorem. ∎

## 3. Finite sum of the product of Hankel operators

In this section we will discuss when the finite sum of the products of Hankel operators with scalar symbols is zero. This is needed in the next section for characterizing $F$ and $G$ in $L_{n\times n}^{\infty}$ such that $T_{FG} = T_F T_G$.

Let $M_{n\times n}$ be the set of $n \times n$ matrices. Let $A = [a_{ij}] \in M_{n\times n}$, define

$$\|A\|_{\infty} = sup_{1\leq i,j\leq n}|a_{ij}|$$

and $(M_{n\times n})_1$ denotes the closed unit ball of $M_{n\times n}$ in the above norm. Let $P_n$ be the set of $n \times n$ permutation matrices.

**Proposition 4.** *Let* $f_k = (f_{k1}, \cdots, f_{kn})^T$ *for* $k = 1, \cdots, m$ *and* $g = (g_1, \cdots, g_n)^T \in L^2(C^n)$. *Let*

$$S_{kn} := f_{k1} \otimes g_1 + \cdots + f_{kn} \otimes g_n, \quad k = 1, \cdots, m.$$

*Then* $S_{kn} = 0$ *for all* $k = 1, \cdots, m$ *if and only if there are a matrix* $A \in (M_{n\times n})_1$ *and a permutation matrix* $R$ *such that*

$$(R-A)f_k = 0, k = 1, \cdots, m \quad and \quad A^*g = 0.$$

*Proof.* We first prove that $S_{kn} = 0, k = 1, \cdots, m$ if and only if there are a matrix $A_0 \in (M_{n\times n})_1$ and a permutation $\sigma$ such that

$$(I - A_0)f_{k\sigma} = 0, k = 1, \cdots, m \quad and \quad A_0^*g_{\sigma} = 0,$$

where $f_{k\sigma} = (f_{k\sigma(1)}, \cdots, f_{k\sigma(n)})^T$. For any $A \in (M_{n\times n})_1$, any permutation $\sigma$, set

$$x_k := (x_{k1}, \cdots, x_{kn})^T = (I - A)f_{k\sigma}, k = 1, \cdots, m \text{ and } y := (y_1, \cdots, y_n)^T = A^*g_{\sigma},$$

then we have

$$S_{kn} = \sum f_{k\sigma_i} \otimes g_{\sigma_i} = x_{k1} \otimes g_{\sigma(1)} + \cdots + x_{kn} \otimes g_{\sigma(n)} + f_{k\sigma(1)} \otimes y_1 + \cdots + f_{k\sigma(n)} \otimes y_n.$$

The sufficiency follows from the above relation. To prove the necessity we use induction. It is clear that for $n = 1$, the result is true with $A = 1$ or $A = 0$. Now assume the result is true for $n - 1$. Without loss of generality, assume that

$$max_{1\leq i\leq n}\|g_i\|_2 = \|g_j\|_2 > 0$$

for some $j$. Note that if $S_{kn} = 0, k = 1, \cdots, m$, then

$$S_{kn}g_j = \sum_{i=1}^{n} <g_j, g_i> f_i = 0, k = 1, \cdots, m.$$

That is

$$f_{kj} + \sum_{i\neq j} a_j f_{ki} = 0,$$

where $a_i = < g_i, g_j > / < g_i, g_j >, |a_i| \leq 1$ for $i \neq j$. Now we rewrite $S_{kn}$ as

$$S_{kn} = (f_{kj} + \sum_{i \neq j} a_j f_{ki}) \otimes g_j + \sum_{i \neq j} f_{ki} \otimes (g_i - \overline{a_j} g_j).$$

¿From above analysis we have that

$$\sum_{i \neq j} f_{ki} \otimes (g_i - \overline{a_j} g_j) = 0, k = 1, \cdots, m.$$

By induction there exist $A_1 \in (M_{n-1 \times n-1})_1$ and a permutation $\omega$ of $\{1, \cdots, j-1, j+1, \cdots, n\}$ such that

$$(I - A_1)(f_{k\omega(1)}, \cdots, f_{k\omega(j-1)}, f_{k\omega(j+1)}, \cdots, f_{k\omega(n)})^T = 0, k = 1, \cdots, m,$$

$$A_1^*(g_{\omega(1)} - \overline{a_{\omega(1)}} g_{\omega(j)}, \cdots, g_{\omega(j-1)} - \overline{a_{\omega(j-1)}} g_{\omega(j)},$$
$$g_{\omega(j+1)} - \overline{a_{\omega(j+1)}} g_{\omega(j)}, \cdots, g_{\omega(n)} - \overline{a_{\omega(n)}} g_{\omega(j)})^T = 0.$$

Let

$$A_0 = \begin{bmatrix} 0 & a \\ 0 & A_1 \end{bmatrix}, \quad \text{where} \quad a = [-a_{\omega(1)}, \cdots, -a_{\omega(j-1)}, -a_{\omega(j+1)}, \cdots, -a_{\omega(n)}].$$

Take $\sigma$ to be such that $\sigma(1) = j, \sigma(i) = \omega(i+1)$ for $2 \leq i \leq j-1$ and $\sigma(i) = \omega(i)$ for $j+1 \leq i \leq n$. It is easy to check that such $A_0$ and $\sigma$ are what we need. Now let $R$ be the permuation matrix such that $f_{k\sigma} = R f_k$ for all $k = 1, \cdots, m$ and $A = R^* A_0$, then it is easy to check that for such $A$ and $R$

$$(R - A)f_k = 0, k = 1, \cdots, m \quad \text{and} \quad A^* g = 0$$

if and only if

$$(I - A_0)f_{k\sigma} = 0, k = 1, \cdots, m \quad \text{and} \quad A_0^* g_\sigma = 0.$$

The proof is complete. ∎

Next we discuss when the finite sum of the products of Hankel operators is zero.

**Theorem 5.** *Let $f_k = (f_{k1}, \cdots, f_{kn})^T, k = 1, \cdots, m$ and $g = (g_1, \cdots, g_n)^T$ for $f_{ij}$ and $g_i$ in $L^\infty$. Let*

$$T_{kn} := \sum_{i=1}^n H_{f_{ki}}^* H_{g_i}, \quad k = 1, \cdots, m.$$

*Then $T_{kn} = 0$ for all $k = 1, \cdots, m$ if and only if there are a matrix $A \in (M_{n \times n})_1$ and a permutation matrix $R$ such that*

$$(1) \qquad (R - A)f_k \in H_{n \times 1}^\infty, k = 1, \cdots, m \quad \text{and} \quad A^* g \in H_{n \times 1}^\infty.$$

*Proof.* To prove the necessity, we recall the following identity proved in Lemma 1.

$$(2) \qquad V^{-1}[\sum_{i=1}^{n}(H_{f_{ki}}^{*}H_{g_i} - T_{\phi_z}^{*}H_{f_{ki}}H_{g_i}T_{\phi_z})]V = \sum_{i=1}^{n}H_{f_{ki}}k_z \otimes H_{g_i}k_z, k = 1, \cdots, m.$$

Therefore if $\sum_{i=1}^{n}H_{f_{ki}}^{*}H_{g_i} = 0$, then

$$\sum_{i=1}^{n}H_{f_{ki}}k_z \otimes H_{g_i}k_z = 0.$$

In particular for $z = 0$ (i.e., $k_0 = 1$), by Proposition 4, the above equation implies that there exist a matrix $A \in (M_{n \times n})_1$ and a permutation matrix $R$ such that

$$(R - A)[H_{f_{k1}}1, \cdots, H_{f_{kn}}1]^{T} = 0, k = 1, \cdots, m \text{ and } A^{*}[H_{g_1}1, \cdots, H_{g_n}1]^{T} = 0.$$

That is (1) holds. To prove the sufficiency, as in the proof of Proposition 4, we note that for any $A \in (M_{n \times n})_1$, any permutation matrix $R$, if we set

$$x_k = (x_{k1}, \cdots, x_{kn})^{T} = (R - A)f_{k\sigma} \text{ and } y = (y_1, \cdots, y_n)^{T} = A^{*}g_{\sigma},$$

then we have

$$(3) \qquad T_{kn} = H_{x_{k1}}^{*}H_{g_1} + \cdots + H_{x_{kn}}^{*}H_{g_n} + H_{f_{k1}}^{*}H_{y_1} + \cdots + H_{f_{kn}}^{*}H_{y_n}, k = 1, \cdots, m.$$

The above formula and the fact that $H_b$ is zero when $b \in H^{\infty}$ prove the sufficiency part of our theorem. ∎

## 4. Zero semi-commutator or commutator

Brown and Halmos [4] showed that the semi-commutator $T_{\phi\psi} - T_{\phi}T_{\psi}$ is zero exactly when either $\phi$ or $\psi$ is analytic for the scalar functions $\phi$ and $\psi$. Halmos [12] also characterized when the commutator $T_{\phi}T_{\psi} - T_{\psi}T_{\phi}$ is zero. In this section, we will characterize when the semi-commutator $T_{FG} - T_F T_G$ or the commutator $T_F T_G - T_G T_F$ is zero for block Toeplitz operators with matrix symbols $F$ and $G$.

Let $E_j$ be the $n \times n$ matrix unit with $(j, j)$-th entry equal to one and all other entries equal to zero. Note that for a $m \times n$ matrix $B$, $BE_j$ is basically the j-th column of $B$.

**Theorem 6.** *Let $F, G \in L_{n \times n}^{\infty}$. The following are equivalent.*
*(1) The semi-commutator $T_F T_G - T_{FG} (= -H_{F^*}^{*}H_G)$ is zero.*
*(2) There exist matrices $A_j \in (M_{n \times n})_1$ and $R_j \in P_n, j = 1, \cdots, n$ such that*

$$(R_j - A_j)F^{*} \in H_{n \times n}^{\infty} \quad and \quad A_j^{*}GE_j \in H_{n \times n}^{\infty}, j = 1, \cdots, n.$$

*(3)*

$$[|(F_+)^{*} - (F_+)^{*}(z)|^2(z)]^{1/2}[|G_- - G_-(z)|^2(z)]^{1/2} = 0$$

*for all $z \in D$.*
*(4)*

$$[|(F_+)^{*} - (F_+)^{*}(z)|^2(z)]^{1/2}[|G_- - G_-(z)|^2(z)]^{1/2} = 0$$

*for some $z \in D$.*

*Proof.* (1)$\Longleftrightarrow$(2). Let $F = (f_{ij})$ and $G = (g_{ij})$. Since $T_F T_G - T_{FG} = H_{F^*}^* H_G$, $T_F T_G - T_{FG} = 0$ if and only if for each $j = 1, \cdots, n$,

$$\sum_{i=1}^{n} H_{\tilde{f}_{ik}}^* H_{g_{ij}} = 0, \quad k = 1, \cdots, n.$$

By Theorem 5, this is equivalent to that for each $j = 1, \cdots, n$, there exist matrix $A_j \in (M_{n \times n})_1$ and $R_j \in P_n$ such that

$$(R_j - A_j) F^* \in H_{n \times n}^{\infty} \quad \text{and} \quad A_j^* G E_j \in H_{n \times n}^{\infty}.$$

(1)$\Longrightarrow$(3). Since $T_F T_G - T_{FG} = H_{F^*}^* H_G$, by Lemma 2, for all $z \in D$, we have

$$trace\{[H_{F^*}^* H_G - T_{\Phi_z}^*[H_{F^*}^* H_G] T_{\Phi_z}]^*[H_{F^*}^* H_G - T_{\Phi_z}^*[H_{F^*}^* H_G] T_{\Phi_z}]\}$$
$$= tr[(|(F_+)^* - (F_+)^*(z)|^2(z))(|G_- - G_-(z)|^2(z))]$$
$$= tr[[|(F_+)^* - (F_+)^*(z)|^2(z)]^{1/2}[|G_- - G_-(z)|^2(z)]^{1/2} \times$$
$$[|(F_+)^* - (F_+)^*(z)|^2(z)]^{1/2}[|G_- - G_-(z)|^2(z)]^{1/2})^*].$$

Hence

$$[|(F_+)^* - (F_+)^*(z)|^2(z)]^{1/2}[|G_- - G_-(z)|^2(z)]^{1/2} = 0.$$

(3)$\Longrightarrow$(4). This is obvious.

(4)$\Longrightarrow$(1). For a given $z \in D$, define a unitary operator $U_z$ on $H^2$ by $U_z h = h \circ \phi_z k_z$. Let $\mathcal{U}_z = diag\{U_z, \cdots, U_z\}$. Then it is easy to check that

$$\mathcal{U}_z^* T_F \mathcal{U}_z = T_{F \circ \phi_z}.$$

Therefore

$$\mathcal{U}_z^*[H_{F^*}^* H_G - T_{\Phi_z}^* H_{F^*}^* H_G T_{\Phi_z}] \mathcal{U}_z$$
$$= [H_{F^* \circ \phi_z}^* H_{G \circ \phi_z} - T_{\Phi_0}^* H_{F^* \circ \phi_z}^* H_{G \circ \phi_z} T_{\Phi_0}].$$

So it is sufficient to prove that $H_{F^*}^* H_G = 0$ if we assume that

$$[|(F_+)^* - (F_+)^*(0)|^2(0)]^{1/2}[|G_- - G_-(0)|^2(0)]^{1/2} = 0.$$

By Lemma 2, we have

$$trace\{[H_{F^*}^* H_G - T_{\Phi_0}^* H_{F^*}^* H_G T_{\Phi_0}]^*[H_{F^*}^* H_G - T_{\Phi_0}^* H_{F^*}^* H_G T_{\Phi_0}]\}$$
$$= tr[[|(F_+)^* - (F_+)^*(0)|^2(0)]^{1/2}[|G_- - G_-(0)|^2(0)]^{1/2} \times$$
$$([|(F_+)^* - (F_+)^*(0)|^2(0)]^{1/2}[|G_- - G_-(0)|^2(0)]^{1/2})^*].$$

So

$$H_{F^*}^* H_G - T_{\Phi_0}^* H_{F^*}^* H_G T_{\Phi_0} = 0.$$

Thus it follows from a theorem [6] that there is a matrix valued function $M$ in $L_{n \times n}^{\infty}$ such that $H_{F^*}^* H_G = T_M$. But by the Douglas theorem [7] we have that $M = 0$. Hence $H_{F^*}^* H_G = 0$. ∎

Next we study the commutator $T_F T_G - T_G T_F$ by reducing it to the semi-commutator case. To see this, first note that

$$T_F T_G - T_G T_F = T_F T_G - T_{FG} + T_{GF} - T_G T_F + T_{(FG-GF)}$$

$$= -(H_{F^*}^* H_G - H_{G^*}^* H_F) + T_{(FG-GF)}.$$

Let

$$B = \begin{bmatrix} F & -G \\ 0 & 0 \end{bmatrix}, \quad C = \begin{bmatrix} G & 0 \\ F & 0 \end{bmatrix}.$$

A simple calculation gives that

$$H_{B^*}^* H_C = \begin{bmatrix} H_{F^*}^* H_G - H_{G^*}^* H_F & 0 \\ 0 & 0 \end{bmatrix}.$$

Therefore, by the Douglas theorem [7], $T_F T_G - T_G T_F = 0$ if and only if $H_{B^*}^* H_C = 0$ and $FG = GF$. A straightforward computation shows that

$$\begin{aligned}
&(4) \quad |(B_+)^* - (B_+)^*(z)|^2(z) \\
&= \begin{bmatrix} |(F_+)^* - (F_+)^*(z)|^2(z) & -(((F_+)^* - (F_+)^*(z))(G_+ - G_+(z)))(z) \\ -(((G_+)^* - (G_+)^*(z))(F_+ - F_+(z)))(z) & |(G_+)^* - (G_+)^*(z)|^2(z) \end{bmatrix}
\end{aligned}$$

$$\begin{aligned}
&(5) \quad |C_- - C_-(z)|^2(z) \\
&= \begin{bmatrix} |G_- - G_-(z)|^2(z) & ((G_- - G_-(z))((F_-)^* - (F_-)^*(z)))(z) \\ ((F_- - F_-(z))((G_-)^* - (G_-)^*(z)))(z) & |F_- - F_-(z)|^2(z) \end{bmatrix}
\end{aligned}$$

This leads to the following result.

**Theorem 7.** *Let $F, G \in L_{n \times n}^\infty$. The following are equivalent.*
*(1) The commutator $T_F T_G - T_G T_F$ is zero.*
*(2) $GF = FG$ and there exist matrices $A_j \in (M_{2n \times 2n})_1$ and $R_j \in P_{2n}, j = 1, \cdots, n$ such that*

$$(R_j - A_j) \begin{bmatrix} F^* \\ -G^* \end{bmatrix} \in H_{2n \times n}^\infty \quad and \quad A_j^* \begin{bmatrix} G \\ F \end{bmatrix} E_j \in H_{2n \times n}^\infty, \quad j = 1, \cdots, n.$$

*(3) $GF = FG$. And*

$$\begin{bmatrix} |(F_+)^* - (F_+)^*(z)|^2(z) & -(((F_+)^* - (F_+)^*(z))(G_+ - G_+(z)))(z) \\ -(((G_+)^* - (G_+)^*(z))(F_+ - F_+(z)))(z) & |(G_+)^* - (G_+)^*(z)|^2(z) \end{bmatrix}^{1/2} \times$$

$$\begin{bmatrix} |G_- - G_-(z)|^2(z) & ((G_- - G_-(z))((F_-)^* - (F_-)^*(z)))(z) \\ ((F_- - F_-(z))((G_-)^* - (G_-)^*(z)))(z) & |F_- - F_-(z)|^2(z) \end{bmatrix}^{1/2} = 0$$

*for some $z \in D$.*

We remark that as in Theorem 6, the matrix in the statement (3) of Theorem 7 is zero for some $z \in D$ if and only if it is zero for all $z \in D$.

An operator $A$ is said to be normal if $AA^* - A^*A = 0$. We observe that by taking $G = F^*$ in Theorem 7 and noting that $T_F^* = T_{F^*}$, we have the following characterization of normal block Toeplitz operators.

**Corollary 8.** *Let $F \in L_{n \times n}^\infty$. The following are equivalent.*
*(1) $T_F$ is normal.*
*(2) $F^*F = FF^*$. And there exist matrices $A_j \in (M_{2n \times 2n})_1$ and $R_j \in P_{2n}, j = 1, \cdots, n$ such that*

$$(R_j - A_j) \begin{bmatrix} F^* \\ -F \end{bmatrix} \in H_{2n \times n}^\infty \quad \text{and} \quad A_j^* \begin{bmatrix} F^* \\ F \end{bmatrix} E_j \in H_{2n \times n}^\infty, \quad j = 1, \cdots, n.$$

*(3) $F^*F = FF^*$. And*

$$\begin{bmatrix} |(F_+)^* - (F_+)^*(z)|^2(z) & -((F_- - F_-(z))(F_+ - F_+(z))(z) \\ -(((F_+)^* - (F_+)^*(z))((F_-)^* - (F_-)^*(z))(z) & |F_- - F_-(z)|^2(z) \end{bmatrix}^{1/2} \times$$

$$\begin{bmatrix} |(F_+)^* - (F_+)^*(z)|^2(z) & ((F_- - F_-(z))(F_+ - F_+(z))(z) \\ (((F_+)^* - (F_+)^*(z))((F_-)^* - (F_-)^*(z))(z) & |F_- - F_-(z)|^2(z) \end{bmatrix}^{1/2} = 0$$

*for some $z \in D$.*

*Proof.* Let $|(B_+)^* - (B_+)^*(z)|^2(z)$ and $|C_- - C_-(z)|^2$ be defined by (4) and (5) with $G = F^*$. The corollary follows from Theorem 7 by noting that

$$(G_+)^* - (G_+)^*(z) = F_- - F_-(z), \quad G_- - G_-(z) = (F_+)^* - (F_+)^*(z). \qquad \blacksquare$$

## 5. A distribution function inequality

Recall that a necessary condition for the compactness of the semi-commutator $T_{FG} - T_F T_G$ is obtained in Theorem 3. Namely, the compactness of $T_{FG} - T_F T_G$ implies that

$$\lim_{|z| \to 1} \|[[|(F_+)^* - (F_+)^*(z)|^2(z)]^{1/2}[|G_- - G_-(z)|^2(z)]^{1/2}\| = 0.$$

To prove that the above condition is also a sufficient condition for the compactness of $T_{FG} - T_F T_G$, we need a certain distribution function inequality. The distribution function inequality involves the Lusin area integral and the Hardy-Littlewood maximal function. The idea to use distribution function inequalities in the theory of Toeplitz operators and Hankel operators first appeared in [2]. In this section we will get such a distribution function inequality.

For $w$ a point of $\partial D$, we let $\Gamma_w$ denote the angle with vertex $w$ and opening $\pi/2$ which is bisected by the radius to $w$. The set of points z in $\Gamma_w$ satisfying $|z - w| < \epsilon$ will be denoted by $\Gamma_{w,\epsilon}$. For h in $L^1(\partial D)$, we define the Lusin area integral of h to be

$$A_\epsilon(h)(w) = [\int_{\Gamma_{w,\epsilon}} |grad h(z)|^2 dA(z)]^{1/2}$$

where h(z) means the harmonic extension of h on D and $dA(z)$ denotes the Lebesgue measure on the unit disk $D$. The Hardy-Littlewood maximal function of the function h will be denoted $Mh$, and for $r > 1$, we let $\Lambda_r h = [M|h|^r]^{1/r}$. For $z \in D$, we let $I_z$ denote the closed subarc of $\partial D$ with center $\frac{z}{|z|}$ and measure $\delta(z) = 1 - |z|$. The Lebesgue measure of the subset E of $\partial D$ will be denoted by $|E|$.

Let $f_i$ and $g_i$ $(i = 1, \cdots, n)$ be in $L^\infty(\partial D)$ and u and v be in $H^2(\partial D)$. We define a generalized area integral to be

$$B_\epsilon(u, v)(w) = \int_{\Gamma_{w,\epsilon}} |\sum_{i=1}^n grad(H_{f_i}u)\overline{grad(H_{g_i}v)}|dA(z).$$

For $l > 2$ and $z \in D$, define

$$\Xi_l(z) = \inf\{\sum_{i=1}^n (\|x_i \circ \phi_z - P(x_i \circ \phi_z)\|_l + \|y_i \circ \phi_z - P(g_i \circ \phi_z)\|_l) : A \in (M_{n \times n})_1, \ R \in P_n\},$$

where

$$x = (x_1, \cdots, x_n)^T = (R - A)f \text{ and } y = (y_1, \cdots, y_n)^T = A^*g,$$

$$\Gamma_l(z) = \sum_{i=1}^n (\|f_i \circ \phi_z - P(f_i \circ \phi_z)\|_l + \|g_i \circ \phi_z - P(g_i \circ \phi_z)\|_l).$$

We have the following distribution function inequality.

**Theorem 9.** *Let $f_i$ and $g_i$ be in $L^\infty$, $i = 1, \cdots, n$, and u and v in the Hardy space $H^2$. Let $z$ be a point in D such that $|z| > 1/2$. Then for any $l > 2$, for $a > 0$ sufficiently large and $\delta(z) = 1 - |z|$*

$$|\{w \in I_z : B_{2\delta(z)}(u, v)(w) < a\Xi_l(z)\Gamma_l(z) \inf_{w \in I_z} \Lambda_r(u)(w) \inf_{w \in I_z} \Lambda_r(v)(w)\}| \geq C_a|I_z|$$

*where $C_a$ depends only on l and a, $lim_{a \to \infty} C_a = 1$, and $1/l + 1/r = 2/p$ for some $1 < p < 2$ and $1 < r < 2$.*

*Proof.* Let $f_i$ and $g_i$ be in $L^\infty$ and u and v in $H^2$. By our definition

$$B_\epsilon(u, v)(w) = \int_{\Gamma_{w,\epsilon}} |\sum_{i=1}^n grad(H_{f_i}u)\overline{grad(H_{g_i}v)}|dA(w).$$

For a fixed $z \in D$, let $\chi_{2I_z}$ denotes the characteristic function of the subset $2I_z$ of $\partial D$ and write $H_{f_i}u = u_{i1} + u_{i2}$, $H_{g_i}v = v_{i1} + v_{i2}$, where

$$u_{i1} = (I - P)[(f_i - P(f \circ \phi_z)) \circ \phi_z)\chi_{2I_z}u],$$

$$u_{i2} = (I - P)[(f_i - P(f_i \circ \phi_z)) \circ \phi_z)(1 - \chi_{2I_z})u],$$

$$v_{i1} = (I - P)[(g_i - P(g_i \circ \phi_z)) \circ \phi_z \chi_{2I_z} v],$$
$$v_{i2} = (I - P)[(g_i - P(g_i \circ \phi_z)) \circ \phi_z (1 - \chi_{2I_z}) v].$$

Thus we have

$$B_\epsilon(u, v)(w) \leq \sum_{i=1}^n (A_\epsilon(u_{i1}) + A_\epsilon(u_{i2}))(A_\epsilon(v_{i1}) + A_\epsilon(v_{i2})) \tag{6}$$

¿From now on we use $\epsilon$ to denote $2\delta(z)$. Let $i$ ($1 \leq i \leq n$) be fixed. We first recall the following estimates from [15] for the terms $A_\epsilon(u_{i1})$ and $A_\epsilon(u_{i2})$ (similarly $A_\epsilon(v_{i1})$ and $A_\epsilon(v_{i2})$).

For $l > 2$, there are a positive constant $C$ and $r \in (1, 2)$ such that

$$[\int_{I_z} A_\epsilon(u_{i1})^p d\sigma(w)]^{1/p} \leq C|I_z|^{1/p} \|f_i \circ \phi_z - P(f_i \circ \phi_z)\|_l \inf_{w \in I_z} \Lambda_r u(w) \tag{7}$$

for some $p > 1$ so that $1/l + 1/r = 1/p$.

For $l > 2$, on $I_z$,

$$A_\epsilon(u_{i2}) \leq C \|f_i \circ \phi_z - P(f_i \circ \phi_z)\|_l \inf_{w \in I_z} \Lambda_{l'} u(w), \tag{8}$$

for some $C > 0$ and $1/l + 1/l' = 1$.

For completeness we give a proof here. We first write $u_{ij} = (I - P)u_j$ ( $j = 1, 2$) where

$$u_1 = (f_{i-} - f_{i-}(z)) \chi_{2I_z} u, \quad u_2 = (f_{i-} - f_{i-}(z))(1 - \chi_{2I_z} u)$$

and $f_{i-} = (I - P)f_i$.

We first prove (7). Note that for $l > 2$, we can always find $l' > 2$ and $p > 1$ so that $l = l'p$ and $r = p \frac{l'}{l'-1} < 2$. By the theorem of Marcinkiewicz and Zygmund, the truncated Lusin area integral $A_\epsilon f(w)$ is $L^p$-bounded for $1 < p < \infty$. So for $l > 2$, we have

$$\int_{I_z} [A_\epsilon(I - P)(u_1)(w)]^p d\sigma(w) \leq C \int_{\partial D} |u_1|^p d\sigma(w)$$

$$= C \int_{2I_z} |f_{i-} - f_{i-}(z)|^p |u(w)|^p d\sigma(w)$$

$$\leq |2I_z| [\frac{1}{|2I_z|} \int_{2I_z} |f_{i-}(w) - f_{i-}(z)|^{pl'} d\sigma(w)]^{1/l'} [\frac{1}{|2I_z|} \int_{2I_z} |u|^r d\sigma(w)]^{p/r}.$$

Let $P(z, w)$ denote the Poisson kernel for the point z. Since

$$[\frac{1}{|2I_z|} \int_{2I_z} |u|^r d\sigma(w)]^{1/r} \leq \Lambda_r u(w)$$

for each $w \in 2I_z$, and an elementary estimate shows that for $w \in 2I_z$, $P(z, w) > \frac{C}{|2I_z|}$, it follows that

$$[\int_{I_z} A_\epsilon((I - P)(u_1)(w))^p d\sigma(w)]^{1/p} \leq C|I_z|^{1/p} [|f_{i-} - f_{i-}(z)|^l(z)]^{1/l} \inf_{w \in I_z} \Lambda_r u(w).$$

Thus (7) follows from the following identity.

$$(9) \qquad [|f_{i-} - f_{i-}(z)|^l(z)]^{1/l} = \|f_i \circ \phi_z - P(f_i \circ \phi_z)\|_l$$

Now we prove (8). For $u_2$, we shall use a pointwise estimate of the norm of the gradient of $(I-P)u_2$. It is easy to see that

$$(I-P)(u_2)(w) = \frac{1}{2\pi} \int \frac{\overline{w}\xi u_2(\xi)}{1 - \overline{w}\xi} d\sigma(\xi).$$

Thus

$$|grad(I-P)u_2(w)| \le C \int \frac{|u_2(\xi)|}{|1 - \overline{w}\xi|^2} d\sigma(\xi)$$

$$\le C \int_{\partial D/2I_z} \frac{|[f_{i-}(\xi) - f_{i-}(z)]u(\xi)|}{|1 - \overline{w}\xi|^2} d\sigma(\xi)$$

On the other hand, there is a constant $C > 0$ so that

$$|\frac{1 - \overline{z}\xi}{1 - \overline{w}\xi}| \ge C$$

for all $\xi$ in $\partial D/2I_z$ and $w$ in $I_z$. Thus we obtain

$$|grad(I-P)u_2(w)| \le C \int_{\partial D/2I_z} \frac{|[f_{i-}(\xi) - f_{i-}(z)]u(\xi)|}{|1 - \overline{z}\xi|^2} d\sigma(\xi).$$

Applying the Hölder inequality yields

$$|grad(I-P)u_2(w)| \le \frac{C}{1 - |z|^2} [|f_{i-} - f_{i-}(z)|^l(z)]^{1/l} [(|u|^{l'})(z)]^{1/l'}.$$

Because the nontangential maximal function is bounded by a constant times the Hardy-Littlewood maximal function, and because z belongs to $\Gamma_{w,\epsilon}$, the last factor on the right is no larger than $C\Lambda_{l'}u(w)$, and again the desired inequality is established by noting (9).

Now we can estimate the products $A_\epsilon(u_{ij})A_\epsilon(v_{ij})$, $j = 1, 2$ by using (7) and (8). By Hölder inequality, we have

$$[\int_{I_z} [A_\epsilon(u_{ij})A_\epsilon(v_{ij})]^{p/2} d\sigma(w)] \le [\int_{I_z} [A_\epsilon(u_{ij})]^p d\sigma(w)]^{1/2} [\int_{I_z} [A_\epsilon(v_{ij})]^p d\sigma(w)]^{1/2}.$$

By using estimate (7) for integral of $A_\epsilon(u_{i1})$ or $A_\epsilon(v_{i1})$ and estimate (8) for integral of $A_\epsilon(u_{i2})$ or $A_\epsilon(v_{i2})$, we get that for $l > 2$, there is $r$ with $1 < r < 2$ such that

$$[\int_{I_z} [A_\epsilon(u_{ij})A_\epsilon(v_{ij})]^{p/2} d\sigma(w)]^{2/p} \le C|I_z|^{2/p} \|f_i \circ \phi_z - P(f_i \circ \phi_z)\|_l \|g_i \circ \phi_z - P(g_i \circ \phi_z)\|_l \times$$

$$(10) \qquad \inf_{w \in I_z} \Lambda_r u(w) \inf_{w \in I_z} \Lambda_r v(w)$$

for some $p > 1$ so that $1/l + 1/r = 2/p$, and a constant C depends on only $l$.

Now we are going to finish the proof by summing up the estimates as above. ¿From (6) we have

$$[\int_{I_z} B(u,v)^{p/2}d\sigma(w)]^{2/p} \leq C(\int_{I_z} \sum_{i=1}^{n}\sum_{j=1}^{2}\{[A_\epsilon(u_{ij})(w)A_\epsilon(v_{ij})(w)]^{p/2})^{2/p}.$$

By estimate (10), we have

$$\int_{I_z} B(u,v)^{p/2}d\sigma(w) \leq C|I_z| \times$$

$$(11) \quad (\sum_{i=1}^{n}\|(f_i\circ\phi_z - P(f_i\circ\phi_z)\|_l^{p/2}\|g_i\circ\phi_z - P(g_i\circ\phi_z)\|_l^{p/2})(\inf_{w\in I_z}\Lambda_s u(w)\inf_{w\in I_z}\Lambda_s v(w))^{p/2}.$$

Next for any $A \in (M_{n\times n})_1$ and any permutation matrix $R \in P_n$, let

$$x = (x_1,\cdots,x_n)^T = (R-A)f \text{ and } y = (y_1,\cdots,y_n)^T = A^*g,$$

where $f = (f_1,\cdots,f_n)^T$ and $g = (g_1,\cdots,g_n)^T$. We apply the above inequality (11) to the vector-valued functions $x$ and $y$. We note that the $B_\epsilon(u,v)$ corresponding to $f$ and $g$ is the same as the $B_\epsilon(u,v)$ corresponding to $x$ and $y$; more precisely,

$$B_\epsilon(u,v)(w) = \int_{\Gamma_{w,\epsilon}} |\sum_{i=1}^{n} grad(H_{f_i}u)\overline{grad(H_{g_i}v)}|dA(w)$$

$$= \int_{\Gamma_{w,\epsilon}} |\sum_{i=1}^{n}(grad(H_{x_i}u)\overline{grad(H_{g_i}v)} + grad(H_{f_i}u)\overline{grad(H_{y_i}v)})|dA(w).$$

By formula (11), we have

$$\int_{I_z} B(u,v)^{p/2}d\sigma(w) \leq C|I_z|(\inf_{w\in I_z}\Lambda_s u(w)\inf_{w\in I_z}\Lambda_s v(w)]^{p/2}) \times$$

$$(\sum_{i=1}^{n}(\|x_i\circ\phi_z - P(x_i\circ\phi_z)\|_l^{p/2}\|g_i\circ\phi_z - P(g_i\circ\phi_z)\|_l^{p/2} +$$

$$\|f_i\circ\phi_z - P(f_i\circ\phi_z)\|_l^{p/2}\|y_i\circ\phi_z - P(y_i\circ\phi_z)\|_l^{p/2})).$$

Therefore

$$[\int_{I_z} B_\epsilon(u,v)^{p/2}d\sigma(w)]^{2/p} \leq C|I_z|^{2/p}\Xi_l(z)\Gamma_l(z)\inf_{w\in I_z}\Lambda_s u(w)\inf_{w\in I_z}\Lambda_s v(w).$$

Next for a fixed z in D and $a > 0$, let $E(a)$ be the set of points in $I_z$ where

$$B_\epsilon(u,v) \leq a\Xi_l(z)\Gamma_l(z)\inf_{w\in I_z}\Lambda_s h(w)\inf_{w\in I_z}\Lambda_s v(w).$$

Then

$$|I_z/E(a)|^{2/p} a \Xi_l(z) \Gamma_l(z) \inf_{w \in I_z} \Lambda_s h(w) \inf_{w \in I_z} \Lambda_s v(w) \leq [\int_{I_z} B_\epsilon(u,v) d\sigma(w)]^{2/p}$$

$$\leq C |I_z|^{2/p} \Xi_l(z) \Gamma_l(z) \inf_{w \in I_z} \Lambda_s h(w) \inf_{w \in I_z} \Lambda_s v(w).$$

So

$$|I_z/E(a)| \leq C a^{-p/2} |I_z|.$$

Therefore for a sufficient large $a > 0$, we have

$$|E(a)| \geq (1 - C a^{-p/2}) |I_z|.$$

Let $C_a = 1 - C a^{-p/2}$. This completes the proof of the theorem. ∎

## 6. Compact finite sum of products

Before proceeding to our main results in this section, we need to introduce some notations involving the maximal ideal space of an algebra. Let $\mathcal{M}$ be the maximal ideal space of $H^\infty$, which is defined to be the set of multiplicative linear maps from $H^\infty$ onto the field of complex numbers. Each multiplicative linear functional $\phi \in \mathcal{M}$ has norm 1 (as an element of the dual of $H^\infty$). If we think of $\mathcal{M}$ has a subset of the dual space $H^\infty$ with weak-star topology then $\mathcal{M}$ becomes a compact Hausdorff space. For $z \in D$ the evaluation functional $f \to f(z)$ is a multiplicative functional. So we can think of $D$ as a subset of $\mathcal{M}$. The Carleson corona theorem tells us that $D$ is dense in $\mathcal{M}$.

By using the Gelfand transform, we can think of $H^\infty$ as a subset of $C(\mathcal{M})$, the continuous, complex-valued functions on the maximal ideal space of $H^\infty$. Explicitly, for $f \in H^\infty$, we extend f from $D$ to $\mathcal{M}$ by defining

$$f(\tau) = \tau(f)$$

for every $\tau \in \mathcal{M}$. Note that this definition is consistent with our earlier identification of $D$ with a subset of $\mathcal{M}$.

By the Hahn-Banach theorem each $\tau \in \mathcal{M}$ extends to a linear functional $\tau'$ on $L^\infty$. In fact, there is a unique representing measure $d\mu$ supported on $M(L^\infty)$, the maximal ideal space of $L^\infty$, such that for each $g \in L^\infty$, $\tau'(g) = \int_{supp(\tau')} g d\mu$. A subset of $M(L^\infty)$ will be called a support set, denoted by $supp\tau$, if it is the (closed) support set of the representing measure for the extension of a functional $\tau$ in $M(H^\infty + C)$.

For $f \in L^\infty$, we let $H^\infty[f]$ denote the closed subalgebra of $L^\infty$ generated by $H^\infty$ and the function $f$. If $f = (f_1, \cdots, f_n)^T$, we still use $H^\infty[f]$ to denote the closed subalgebra of $L^\infty$ generated by $H^\infty$ and functions $f_1, \cdots, f_n$. Recall that

$$\Xi_2(z) = \inf\{\sum_{i=1}^n (\|x_i \circ \phi_z - P(x_i \circ \phi_z)\|_2 + \|y_i \circ \phi_z - P(g_i \circ \phi_z)\|_2) : A \in (M_{n \times n})_1 \ R \in P_n\},$$

where

$$x = (x_1, \cdots, x_n)^T = (R - A)f \text{ and } y = (y_1, \cdots, y_n)^T = A^* g.$$

**Theorem 10.** *Let $f = (f_1, \cdots, f_n)^T$ and $g = (g_1, \cdots, g_n)^T$ for $f_i$ and $g_i$ in $L^\infty$. The following are equivalent.*
*(1) $H^*_{f_1} H_{g_1} + \cdots + H^*_{f_n} H_{g_n}$ is compact.*
*(2) $\lim_{z \to \partial D} \|\sum_{i=1}^n (H_{f_i} k_z) \otimes (H_{g_i} k_z)\| = 0$.*
*(3) $\lim_{z \to \partial D} \Xi_2(z) = 0$.*
*(4) For each $m \in M(H^\infty + C)$, there exist a matrix $A \in (M_{n \times n})_1$ and a permutation matrix $R \in P_n$ such that $(R-A)f|\mathrm{suppm} \in H^\infty|\mathrm{suppm} \quad and \quad A^*g|\mathrm{suppm} \in H^\infty|\mathrm{suppm}$*
*(5) The following relation holds.*

$$\cap_{\{A \in (M_{n \times n})_1, \, R \in P_n\}} H^\infty[(R-A)f, A^*g] \subset H^\infty + C.$$

*Proof.* Without loss of generality we may assume that $\|f_i\|_\infty < 1/2$ and $\|g_i\|_\infty < 1/2$ for all $i = 1, \cdots, n$.

(1) $\Longrightarrow$ (2). Assume that $\sum_{i=1}^n H^*_{f_i} H_{g_i}$ is compact. By Lemma 2 [15] we obtain

$$(12) \qquad \lim_{|z| \to 1} \|\sum_{i=1}^n H^*_{f_i} H_{g_i} - T^*_{\phi_z} (\sum_{i=1}^n H^*_{f_i} H_{g_i}) T_{\phi_z}\| = 0.$$

But by the proof of Lemma 1,

$$V^{-1} (\sum_{i=1}^n H^*_{f_i} H_{g_i} - T^*_{\phi_z} (\sum_{i=1}^n H^*_{f_i} H_{g_i}) T_{\phi_z}) V = \sum_{i=1}^n H_{f_i} k_z \otimes H_{g_i} k_z,$$

where recall that $V$ is the antiunitary operator on $L^2(\partial \mathbb{D})$ defined by $Vh(e^{i\theta}) = e^{-i\theta} \overline{h(e^{i\theta})}$. Thus

$$(13) \qquad \lim_{|z| \to 1} \|\sum_{i=1}^n H_{f_i} k_z \otimes H_g k_z\| = 0.$$

That is (2) holds.

(2) $\Longrightarrow$ (3). Assume now that (2) holds. Suppose that (3) does not hold. That is, there are $\delta > 0$ and a net $\{z\} \subset D$ accumulating a point in $\partial D$ such that

$$\Xi_2(z) \geq \delta.$$

We will get a contradiction. We may assume that the net $\{z\}$ converges to some nontrivial point $m \in M(H^\infty + C)$.

Let $H^\infty_m$ denote the algebra $H^\infty|_{\mathrm{suppm}}$ on $\mathrm{suppm}$, and $L^\infty_m$ denote the algebra $L^\infty|_{\mathrm{suppm}}$. Then $L^\infty_m / H^\infty_m$ is a vector space. For a function $\psi$ in $L^\infty$, let $[\psi]_m$ denote the element in $L^\infty_m / H^\infty_m$ which contains $\psi$. For $f = (f_1, \cdots, f_n)^T$, let $[f]_m = ([f_1]_m, \cdots, [f_n]_m)^T$ and $f \in H^\infty_m$ means that $f_i \in H^\infty_m$ for all $i = 1, \cdots, n$. Let $g = (g_1, \cdots, g_n)^T$. For convenience we also introduce the following notations.

$$\|H_f k_z\|_2 := \sum_{i=1}^n \|H_{f_i} k_z\|_2,$$

$$H_f k_z \otimes H_g k_z := \sum_{i=1}^{n} H_{f_i} k_z \otimes H_{g_i} k_z.$$

Suppose that the dimension of the space spanned by $[f_1]_m, ..., [f_n]_m$ is $N \le n$. We may assume that $\{[f_1]_m, ..., [f_N]_m\}$ is a basis such that $([f_1]_m, ..., [f_n]_m)^T = B([f_1]_m, ..., [f_N]_m)^T$ with $B = (b_{ij})$ and $|b_{ij}| \le 1$, see the proof of Proposition 4 for details. Let $A$ be the matrix $(B, 0)_{n \times n}$. Then $f - Af$ is in $H_m^\infty$ on the support set $supp m$. By Lemma 3 [15],

$$\lim_{z \to m} \|H_{f-Af} k_z\|_2 = 0.$$

On the other hand,

$$H_f k_z \otimes H_g k_z = H_{f-Af} k_z \otimes H_g k_z + H_f k_z \otimes H_{A^*g} k_z$$
$$= H_{f-Af} k_z \otimes H_g k_z + H_{f_{(N)}} k_z \otimes H_{B^*g} k_z,$$

where $f_{(N)} = (f_1, ..., f_N)$. As $z$ goes to $m$, the first term in the right hand side of the above equation goes to zero. Hence

$$\lim_{z \to m} \|H_{f_{(N)}} k_z \otimes H_{B^*g} k_z\| = 0.$$

We are going to show that

$$\lim_{z \to m} \|H_{B^*g} k_z\|_2 = 0.$$

Suppose that this is not true. We may assume that

$$\lim_{z \to m} \|H_{g_1} k_z\|_2 > 0.$$

Let $a_i(z) = \langle H_{g_1} k_z, H_{g_i} k_z \rangle$. Note that $|a_i(z)| \le \sqrt{\|g_1\|_\infty} \sqrt{\|g_i\|_\infty} \le 1$. We may assume that $a_i(z)$ converges to $a_i$ as $z$ goes to $m$. By our assumption $a_1 \ne 0$. But

$$\lim_{z \to m} \|(H_{f_{(N)}} k_z \otimes H_{B^*g} k_z) H_{g_1} k_z\|_2 = 0$$

implies that

$$\lim_{z \to m} \|H_{\sum_{i=1}^{N} a_i f_i} k_z\|_2 = 0.$$

By Lemma 3 [15], $\sum_{i=1}^{N} a_i f_i$ is in $H^\infty$ on $supp m$. This contradicts the fact that $\{[f_1]_m, ..., [f_N]_m\}$ is a basis. Therefore

$$\lim_{z \to m} \|H_{B^*g} k_z\|_2 = 0.$$

Hence

$$\lim_{z \to m} \|H_{A^*g} k_z\|_2 = 0.$$

But

$$\|H_{(I-A)f} k_z\| + \|H_{A^*g} k_z\| \ge \Xi_2(z).$$

Hence

$$\lim_{z \to m} (\|H_{(I-A)f} k_z\| + \|H_{A^*g} k_z\|) \ge \delta.$$

This is a contradiction.

$(3) \Longrightarrow (4)$. We are going to show that for each $m \in M(H^\infty + C)$, there exist matrices $A_m \in (M_{n \times n})_1$ and $R_m \in P_n$ such that

$$[(R_m - A_m)f]_m = 0, \ [A_m^* g]_m = 0.$$

Let $z$ be a net in $D$ converging to $m$. By condition (3), there are matrices $A_z \in (M_{n \times n})_1$ and $R_z \in P_n$ such that

$$\lim_{z \to m} (\|H_{(R_z - A_z)f} k_z\|_2 + \|H_{A_z^* g} k_z\|_2) = 0.$$

Since $(M_{n \times n})_1$ is compact and the permutation group $P_n$ is also compact, we may assume that $A_z$ converges to $A_m$ and $R_z$ converges to $R_m$. Hence

$$\lim_{z \to m} (\|H_{(R_m - A_m)f} k_z\|_2 + \|H_{A_m^* g} k_z\|_2) = 0.$$

By Lemma 3 [15], we have

$$[(R_m - A_m)f]_m = 0, \ [A_m^* g]_m = 0.$$

$(4) \Longrightarrow (5)$. By the Chang-Marshall theorem [8], we need only to show that

$$M(H^\infty + C) \subset M(\cap_{\{A \in (M_{n \times n})_1, \ R \in P_n\}} H^\infty[(R - A)f, A^* g]).$$

Condition (4) states exactly that

$$M(H^\infty + C) \subset \cup_{\{A \in (M_{n \times n})_1, \ R \in P_n\}} M(H^\infty[(R - A)f, A^* g]).$$

By the Sarason theorem [11],

$$M(\cap_{\{A \in (M_{n \times n})_1, \ R \in P_n\}} H^\infty[(R - A)f, A^* g])$$

(14) $$= Closure\ of\ \cup_{\{A \in (M_{n \times n})_1, \ R \in P_n\}} M(H^\infty[(R - A)f, A^* g]).$$

Hence

$$M(H^\infty + C) \subset M(\cap_{\{A \in (M_{n \times n})_1, \ R \in P_n\}} H^\infty[(R - A)f, A^* g]).$$

$(5) \Longrightarrow (3)$. Suppose that (3) does not hold. There are $\delta > 0$ and a net $z$ in $D$ converging to some $m \in M(H^\infty + C)$ such that

$$\Xi_2(z) \geq \delta.$$

By condition (5) and Sarason's Theorem [11] as in (14), there are a net $m_\alpha \in M(H^\infty)$ and matrices $A_\alpha$, $R_\alpha$ such that $m_\alpha$ converges to $m$ and

$$[(R_\alpha - A_\alpha)f]_{m_\alpha} = 0, \ [A_\alpha^* g]_{m_\alpha} = 0.$$

We may assume that $A_\alpha$ converges to some $A_m$ and $R_\alpha$ converges to some $R_m$. We claim that

$$[(R_m - A_m)f]_m = 0, \ [A_m^* g]_m = 0.$$

As in Lemma 3 [15], let $u_i$ $(i=1,\cdots,n)$ be the unimodular functions such that $u_m = (u_1,...,u_n)^T$ in $(R_m - A_m)f + H_n^\infty$ and $u_i$ $(i=1,\cdots,n)$ be the unimodular functions such that $v_m = (v_1,...,v_n)^T$ in $A_m^* g + H_n^\infty$. Then by Lemma 3 [15]

$$\sum_{i=1}^n [(1-|u_i(z)|^2) + (1-|v_i(z)|^2)]$$

$$\leq C[\|A_\alpha - A_m\|_\infty + \|R_\alpha - R_m\|_\infty] + C[\|H_{(R_\alpha - A_\alpha)f} k_z\|_2 + \|H_{A_\alpha^* g} k_z\|_2]$$

for all $z \in D$. Hence

$$\sum_{i=1}^n [(1-|u_i(m_\alpha)|^2) + (1-|v_i(m_\alpha)|^2)] \leq C[\|A_\alpha - A_m\|_\infty + \|R_\alpha - R_m\|_\infty].$$

Since these functions $u_i$ and $v_i$ are continuous on $M(H^\infty)$, we have

$$\sum_{i=1}^n [(1-|u_i(m)|^2) + (1-|v_i(m)|^2)] = 0.$$

Therefore

$$[(R_m - A_m)f]_m = 0, \ \ [A_m^* g]_m = 0.$$

This proves our claim. But again by Lemma 3 [15], this implies that

$$\lim_{z \to m} (\|H_{(R_m - A_m)f} k_z\|_2 + \|H_{A_m^* g} k_z\|_2) = 0.$$

This contradicts to the assumption that

$$\|H_{(R_m - A_m)f} k_z\|_2 + \|H_{A_m^* g} k_z\|_2 \geq \Xi_2(z) > \delta,$$

(3) $\Longrightarrow$ (1). Now we assume that

$$\lim_{z \to \partial D} \Xi_2(z) = 0.$$

We use the distribution inequality obtained in Section 5 to show that $\sum_{i=1}^n H_{f_i}^* H_{g_i}$ is compact. Since the quantity $\Xi_r(z)$ for some $r > 2$ appears in the distribution inequality, we first need to show that in fact for some $r$ such that $3 > r > 2$,

$$\lim_{z \to \partial D} \Xi_r(z) = 0 \tag{15}$$

Recall that

$$\Xi_l(z) = \inf \{ \sum_{i=1}^n (\|x_i \circ \phi_z - P(x_i \circ \phi_z)\|_l + \|y_i \circ \phi_z - P(g_i \circ \phi_z)\|_l) : A \in (M_{n\times n})_1 \ R \in P_n \},$$

where

$$x = (x_1,\cdots,x_n)^T = (R-A)f \text{ and } y = (y_1,\cdots,y_n)^T = A^* g.$$

First note that since $A \in (M_{n \times n})_1$, $R \in P_n$ and by our assumption $\|f_i\|_\infty < 1/2$ and $\|g_i\|_\infty < 1/2$, we have $\|x_i\|_\infty \leq (1+n)/2$ and $\|y_i\|_\infty \leq n/2$. Thus

$$\sum_{i=1}^n (\|x_i \circ \phi_z - P(x_i \circ \phi_z)\|_r + \|y_i \circ \phi_z - P(y_i \circ \phi_z)\|_r)$$

$$\leq C_r (\sum_{i=1}^n (\|x_i \circ \phi_z - P(x_i \circ \phi_z)\|_2 + \|y_i \circ \phi_z - P(y_i \circ \phi_z)\|_2)^{2/r}).$$

for some constant $C_r$ dependent only on $r$ and $n$. Therefore

$$\Xi_r(z) \leq C_r (\Xi_2(z))^{2/r}.$$

Since $\lim_{z \to \partial D} \Xi_2(z) = 0$, we obtain

$$\lim_{z \to \partial D} \Xi_r(z) = 0.$$

This completes the proof of (15).

Now let $u$ and $v$ be two functions in $H^2$. Since $H_{f_i} u$ is orthogonal to $H^2$, we see that $(H_{f_i} u)(0) = 0$. Thus by the Littlewood-Paley formula [8], we have

$$< u, (\sum_{i=1}^n H_{f_i}^* H_{g_i}) v > = \sum_{i=1}^n < H_{f_i} u, H_{g_i} v > =$$

$$\frac{1}{\pi} \int_{\mathbb{D}} (\sum_{i=1}^n grad(H_{f_i} u) \overline{grad(H_{g_i} v)}) \, log \frac{1}{|z|} \, dA(z) = I_R + II_R,$$

where for $1/2 < R < 1$,

$$I_R = \int_{|z| > R} (\sum_{i=1}^n grad(H_{f_i} u) \overline{grad(H_{g_i} v)}) \, log \frac{1}{|z|} \, dA(z)$$

and

$$II_R = \int_{|z| \leq R} (\sum_{i=1}^n grad(H_{f_i} u) \overline{grad(H_{g_i})}) \, log \frac{1}{|z|} \, dA(z).$$

One easily checks that there is a compact operator $K_R$ such that

$$II_R = < u, K_R v > .$$

Thus, if we show that $I_R \to 0$ as $R \to 1$, then $\|T_f T_g - T_g T_f - K_R\| \to 0$, and we are done. The rest of the proof will be devoted to showing that $I_R \to 0$ as $R \to 1$.

Choose $z \in \mathbb{D}$ and fix a constant $a \geq 1$ for which the Distribution Inequality holds; that is

$$|t \in I_z : \{B_\epsilon(u,v)(t) \leq a \Xi_r(z) \Gamma_r(z) \Lambda_s u(t) \Lambda_s v(t)\}| \geq K_a |I_z|.$$

For $t \in \partial \mathbb{D}$, let

$$\rho(t) = max\{\epsilon : B_\epsilon(u, v)(t) \leq a\Xi_r(z)\Gamma_r(z)\Lambda_s u(t)\Lambda_s v(t)\}.$$

Let $\chi_t$ denote the characteristic function of $\Gamma_{t,\rho(t)}$. Then

$$\int_{\partial \mathbb{D}} B_{\rho(t)}(u, v)(t) \, dt \leq a\Xi_r(z)\Gamma_r(z) \int_{\partial \mathbb{D}} \Lambda_s u(t)\Lambda_s v(t) \, dt$$

$$\leq a\Xi_r(z)\Gamma_r(z)||\Lambda_s u||_2||\Lambda_s v||_2.$$

Since $\dfrac{2}{s} > 1$, so by [8]

$$||\Lambda_s u||_2 = ||M(|u|^s)^{1/s}||_2 = [M(|u|^s)||_{2/s}]^{1/s} \leq A_s(||\,|u|^s||_{2/s}|)^{1/s}$$

So

$$||\Lambda_s u||_2 \leq A_s||u||_2.$$

Similarly,

$$||\Lambda_s v||_2 \leq A'_s||v||_2.$$

Thus

(16) $$\int_{\partial \mathbb{D}} B_{\rho(t)}(u, v)(t)dt \leq a^*\Xi_r(z)\Gamma_r(z)||u||_2\,||v||_2.$$

On the other hand,

$$\int_{\partial \mathbb{D}} B_\epsilon(u, v)(t) = \int_{\partial \mathbb{D}} \int_{\Gamma_{t,\rho(t)}} |\sum_{i=1}^n grad(H_{f_i}u)\overline{grad(H_{g_i}v)}|dA(z) \, dt$$

So

$$\int_{\partial \mathbb{D}} B_\epsilon(u, v)(t) \geq \int_{\partial \mathbb{D}} \int_{|z|>R} \chi_\omega(z)|\sum_{i=1}^n grad(H_{f_i}u)\overline{grad(H_{g_i}v)}|dA(z) \, dt.$$

Now the Distribution Function Inequality tells us that $\rho(t) \geq (1-|z|^2)$ on a subset $E_z$ of $I_z$ satisfying

$$|E_z| \geq K_a|I_z|.$$

Now, for $t \in E_z$, we have $t \in I_z$. Thus if we write $z = re^{i\theta}$ and note that $\rho(t) \geq \frac{3}{2}(1-|z|)$ we have

$$|r\,e^{i\theta} - e^{i\,t}| \leq |r\,e^{i\theta} - e^{i\theta}| + |e^{i\theta} - e^{i\,t}| \leq (1-|z|) + \frac{(1-|z|)}{2}) \leq \rho(t)$$

Therefore, for $t \in E_z$, we have that $z \in \Gamma_{t,\rho(t)}$ and that $\chi_t(z) = 1$ on $E_z$. So,

$$\int_{\partial \mathbb{D}} \int_{|z|>R} \chi_t(z)|\sum_{i=1}^n grad(H_{f_i}u)\overline{grad(H_{g_i}v)}|dA(z) \, dt$$

$$\geq \int_{|z|>R} [\int_{\partial \mathbb{D}} \chi_t(z) \, dt]|\sum_{i=1}^n grad(H_{f_i}u)\overline{grad(H_{g_i}v)}|dA(z)$$

Since $\chi_t(z) = 1$ on $E_z$, we have

$$\int_{\partial \mathbb{D}} B_\epsilon(u,v)(t)dt \geq \int_{|z|>R} |E_z| |\sum_{i=1}^n grad(H_{f_i}u)\overline{grad(H_{g_i}v)}|dA(z)$$

But, $|E_z| \geq K_a(1-|z|^2)$, so

$$\int_{\partial \mathbb{D}} B_\epsilon(u,v)(t) \geq K_a \int_{|z|>R} |\sum_{i=1}^n grad(H_{f_i}u)\overline{grad(H_{g_i}v)}|(1-|z|^2)dA(z).$$

Since

$$I_R = \int_{|z|>R} |\sum_{i=1}^n grad(H_{f_i}u)\overline{grad(H_{g_i}v)}| \, \log\frac{1}{|z|} \, dA(z),$$

we have,

$$\int_{\partial \mathbb{D}} B_\epsilon(u,v)(t) \geq K_a|I_R|.$$

Combining this together with (16), we see that $|I_R| \leq C\Xi_r(z)\Gamma_r(z)||u||_2 ||v||_2$. But by (15) ,

$$\lim_{z\to\partial D} \Xi_r(z) = 0.$$

and $\Gamma_r(z)$ is bounded. Hence we have $\lim_{R\to 1} |I_R| = 0$. This finishes the proof. ∎

## 7. Compact semi-commutator or commutator

In this section by combining the results in Sections 2 and 6, we will show several necessary and sufficient conditions for the semi-commutator or the commutator of the block Toeplitz operators with matrix symbols to be compact. We also give a characterization of essentially normal block Toeplitz operators.

**Theorem 11.** *Let $F$ and $G$ be in $L^\infty_{n\times n}$. Let $F^* = (f_1,\cdots,f_n)$ and $G = (g_1,\cdots,g_n)$. The following are equivalent.*
*(1) $T_{FG} - T_F T_G (= H^*_{F^*}H_G)$ is compact.*
*(2) $\cup_{i,j} \cap_{\{A\in(M_{n\times n})_1, \, R\in P_n\}} H^\infty[(R-A)f_i, A^*g_j] \subset H^\infty + C.$*
*(3)*
$$\lim_{|z|\to 1} \|[|(F_+)^* - (F_+)^*(z)|^2(z)]^{1/2}[|G_- - G_-(z)|^2(z)]^{1/2}\| = 0.$$

*Proof.* (1)$\Longleftrightarrow$(2). Let $F = (f_{ij})_{n\times n}$ and $G = (g_{ij})_{n\times n}$. Note that $T_{FG} - T_F T_G$ is compact if and only if each entry of

$$(\sum_k H^*_{\overline{f_{ik}}}H_{g_{kj}})_{ij}$$

is compact. By Theorem 10, this is equivalent to

$$\cup_{i,j} \cap_{\{A\in(M_{n\times n})_1, \, R\in P_n\}} H^\infty[(R-A)f_i, A^*g_j] \subset H^\infty + C.$$

(1)$\Longleftrightarrow$(3). (1)$\Rightarrow$(3) is proved in Theorem 3. We are going to prove that (3)$\Rightarrow$(1). By Lemma 2,

$$tr(|(F_+)^* - (F_+)^*(z)|^2(z))(|G_- - G_-(z)|^2(z)) =$$
$$trace[H_{F^*}^* H_G - T_{\Phi_z}^* H_{F^*}^* H_G T_{\Phi_z}]^*[H_{F^*}^* H_G - T_{\Phi_z}^* H_{F^*}^* H_G T_{\Phi_z}].$$

Thus by Lemma 1, (3) implies that

$$\lim_{|z| \to 1} \| \sum_k H_{\overline{f_{ik}}} k_z \otimes H_{g_{kj}} k_z \| = 0$$

for all $i, j$. By Theorem 10, we have that $\sum_k H_{\overline{f_{ik}}}^* H_{g_{kj}}$ are compact for all $i, j$. Hence $H_{F^*}^* H_G$ is compact. Therefore $T_{FG} - T_F T_G$ is compact. The proof is complete. ∎

Next we characterize when the commutator $T_F T_G - T_G T_F$ is compact. To do this, recall that

$$T_F T_G - T_G T_F == T_F T_G - T_{FG} + T_{GF} - T_G T_F + T_{(FG-GF)}$$

$$= -(H_{F^*}^* H_G - H_{G^*}^* H_F) + T_{(FG-GF)}.$$

Therefore by the Douglas theorem [7], $T_F T_G - T_G T_F$ is compact if and only if $FG = GF$ and $H_{F^*}^* H_G - H_{G^*}^* H_F$ is compact. But if we let

$$B = \begin{bmatrix} F & -G \\ 0 & 0 \end{bmatrix}, \quad C = \begin{bmatrix} G & 0 \\ F & 0 \end{bmatrix},$$

then

$$H_{B^*}^* H_C = \begin{bmatrix} H_{F^*}^* H_G - H_{G^*}^* H_F & 0 \\ 0 & 0 \end{bmatrix}.$$

Therefore, $T_F T_G - T_G T_F$ is compact if and only if $FG = GF$ and $H_{B^*}^* H_C$ is compact. Note that $|(B_+)^* - (B_+)^*(z)|^2(z)$ and $|C_- - C_-(z)|^2(z)$ can be computed as in (4) and (5). The following result now follows immediately from Theorem 11.

**Theorem 12.** *Let $F$ and $G$ be in $L_{n \times n}^\infty$. Let*

$$\begin{bmatrix} F^* \\ -G^* \end{bmatrix} = (f_1, \cdots, f_n), \quad \begin{bmatrix} G \\ F \end{bmatrix} = (g_1, \cdots, g_n).$$

*The following are equivalent.*
*(1) $T_F T_G - T_G T_F$ is compact.*
*(2) $FG = GF$ and $\cup_{i,j} \cap_{\{A \in (M_{n \times n})_1, \, R \in P_n\}} H^\infty[(R - A)f_i, A^* g_j] \subset H^\infty + C.$*
*(3) $FG = GF$ and*

$$\lim_{|z| \to 1} \left\| \begin{bmatrix} |(F_+)^* - (F_+)^*(z)|^2(z) & -(((F_+)^* - (F_+)^*(z))(G_+ - G_+(z)))(z) \\ -(((G_+)^* - (G_+)^*(z))(F_+ - F_+(z)))(z) & |(G_+)^* - (G_+)^*(z)|^2(z) \end{bmatrix}^{1/2} \times \right.$$

$$\left. \begin{bmatrix} |G_- - G_-(z)|^2(z) & ((G_- - G_-(z))((F_-)^* - (F_-)^*(z)))(z) \\ ((F_- - F_-(z))((G_-)^* - (G_-)^*(z)))(z) & |F_- - F_-(z)|^2(z) \end{bmatrix}^{1/2} \right\| = 0.$$

An operator $A$ is said to be essentially normal if $A^*A - AA^*$ is compact. By taking $G = F^*$, we immediately get the following characterization of essentially normal block Toeplitz operators.

**Corollary 13.** *Let $F$ be in $L_{n \times n}^\infty$. Let*

$$\begin{bmatrix} F^* \\ -F \end{bmatrix} = (f_1, \cdots, f_n), \quad \begin{bmatrix} F^* \\ F \end{bmatrix} = (g_1, \cdots, g_n).$$

*The following are equivalent.*
*(1) $T_F$ is essentially normal.*
*(2) $FF^* = F^*F$ and $\cup_{i,j} \cap_{\{A \in (M_{n \times n})_1, R \in P_n\}} H^\infty[(R-A)f_i, A^*g_j] \subset H^\infty + C$.*
*(3) $FF^* = F^*F$ and*

$$\lim_{|z| \to 1} \| \begin{bmatrix} |(F_+)^* - (F_+)^*(z)|^2(z) & -((F_- - F_-(z))(F_+ - F_+(z))(z) \\ -(((F_+)^* - (F_+)^*(z))((F_-)^* - (F_-)^*(z))(z) & |F_- - F_-(z)|^2(z) \end{bmatrix}^{1/2} \times$$

$$\begin{bmatrix} |(F_+)^* - (F_+)^*(z)|^2(z) & ((F_- - F_-(z))(F_+ + F_+(z))(z) \\ (((F_+)^* - (F_+)^*(z))((F_-)^* - (F_-)^*(z))(z) & |F_- - F_-(z)|^2(z) \end{bmatrix}^{1/2} \| = 0.$$

For the scalar symbols there were several other sufficient conditions for the product of two Toeplitz operators to be a compact perturbation of a Toeplitz operator. To state those conditions we need some notations. The fiber $M_\lambda$ of $M(L^\infty)$ above the point $\lambda$ is the set $\{\tau \in M(L^\infty) : z(\tau) = \lambda\}$. We recall that a subset of $M(L^\infty)$ is called an antisymmetric set if any real-valued function in $H^\infty + C$ is constant on the set.

One of the following conditions implies the compactness of the semi-commutator $T_{\bar\phi} T_\psi - T_{\bar\phi\psi}$ of Toeplitz operators with scalar symbols $\phi$ and $\psi$.
**(1)** Either $\phi$ or $\psi$ is in $C(\partial D)$ [5].
**(2)** $\phi$ and $\psi$ are piecewise continuous and have no common discontinuities [10]
**(3)** Either $\overline{\phi}$ or $\psi$ is in $H^\infty$ on each fiber $M_z$ for z on the circle [13].
**(4)** Either $\overline{\phi}$ is in $H^\infty$ or $\psi$ is in $H^\infty$ on each set of maximal antisymmetry of $H^\infty + C$ [1].

It was shown in [2] that $H^\infty[\overline\phi] \cap H^\infty[\psi] \subset H^\infty + C$ is equivalent to
**(5)** Either $\phi$ or $\psi$ is in $H^\infty$ on each support set.

Next we will show some sufficient conditions for the compactness of the semi-commutator $T_F T_G - T_{FG}$ of Toeplitz operators with matrix symbols $F$ and $G$. Those conditions are analogous to the above conditions of the scalar case. Some of them are well known ([7], [9]).

**Corollary 14.** *Let $F$ and $G$ be in $L_{n \times n}^\infty$. Then one of the following conditions is a sufficient condition for $T_F T_G - T_{FG}$ to be compact:*
**(1)** *Either $F^*$ or $G$ is in $C_{n \times n}(\partial D)$.*
**(2)** *$F^*$ and $G$ are piecewise continuous and have no common discontinuities.*
**(3)** *Either $F^*$ or $G$ is in $H_{n \times n}^\infty$ on each fiber $M_z$ for z on the circle.*

**(4)** *Either $F^*$ or $G$ is in $H_{n \times n}^\infty$ on each maximal antisymmetric set of $H^\infty + C$.*

**(5)** *Either $F^*$ or $G$ is in $H_{n \times n}^\infty$ on each support set.*

**(6)** $H_{n \times n}^\infty[F^*] \cap H_{n \times n}^\infty[G] \subset H_{n \times n}^\infty + C_{n \times n}(\partial D)$, *where $H_{n \times n}^\infty[G]$ denotes the subalgebra of $L_{n \times n}^\infty$ generated by $H_{n \times n}^\infty$ and $G$.*

*Proof.* Notice that Conditions (1) to (6) in the corollary are ordered by weakness. So it is sufficient to show that Condition (6) is stronger than Condition (2) in Theorem 11. Since $H_{n \times n}^\infty[G]$ denotes the subalgebra of $L_{n \times n}^\infty$ generated by $H_{n \times n}^\infty$ and $G$, we observe that

$$H_{n \times n}^\infty[G] = (H^\infty[g_{11}, \cdots, g_{1n}, \cdots, g_{n1}, \cdots, g_{nn}])_{n \times n}.$$

Hence Condition (6) is equivalent to

$$H^\infty[\overline{f_{11}}, \cdots, \overline{f_{1n}}, \cdots, \overline{f_{n1}}, \cdots, \overline{f_{nn}}] \cap H^\infty[g_{11}, \cdots, g_{1n}, \cdots, g_{n1}, \cdots, g_{nn}] \subset H^\infty + C.$$

Let $F^* = (f_1, \cdots, f_n)$ and $G = (g_1, \cdots, g_n)$. It is easy to see that

$$\cup_{i,j} \cap_{\{A \in (M_{n \times n})_1, \ R \in P_n\}} H^\infty[(R-A)f_i, A^* g_j]$$

$$\subset \cup_{i,j} H^\infty[f_i] \cap H^\infty[g_j]$$

But for all $i, j$, we have

$$H^\infty[f_i] \cap H^\infty[g_j]$$

$$\subset H^\infty[\overline{f_{11}}, \cdots, \overline{f_{1n}}, \cdots, \overline{f_{n1}}, \cdots, \overline{f_{nn}}] \cap H^\infty[g_{11}, \cdots, g_{1n}, \cdots, g_{n1}, \cdots, g_{nn}]$$

Hence

$$\cup_{i,j} \cap_{\{A \in (M_{n \times n})_1, \ R \in P_n\}} H^\infty[(R-A)f_i, A^* g_j] \subset H^\infty + C.$$

So by Theorem 11, $T_F T_G - T_{FG}$ is compact. This completes the proof. ∎

**Acknowledgment.** We thank S. Treil for his useful suggestions.

CAIXING GU, MATHEMATICS DEPARTMENT, UNIVERSITY OF CALIFORNIA IRVINE, IRVINE, CA 72717
    *E-mail address*: cgu@math.uci.edu

DECHAO ZHENG, MATHEMATICS DEPARTMENT, MICHIGAN STATE UNIVERSITY, EAST LANSING, MI 48824
    *E-mail address*: dechao@math.msu.edu